\documentclass[12pt,a4paper]{article}

\usepackage{amssymb}
\usepackage{amsmath}
\usepackage{indentfirst}
\usepackage{t1enc}
\usepackage{array}
\usepackage[mathscr]{eucal}
\usepackage{enumerate}
\usepackage{graphicx}
\usepackage[colorlinks=true,linkcolor=blue,citecolor=blue,urlcolor=blue]{hyperref}
\usepackage{amsthm}
\usepackage{tikz}

\usepackage{xcolor}
\usepackage[normalem]{ulem}

\usepackage[english]{babel}
\usepackage[utf8]{inputenc}

\theoremstyle{plain}

\newtheorem*{dfn}{\textsc{Definition}}

\begin{document}

\title{The Case of the Second Smallest Modulus: \\ \large or why we cannot be arbitrarily careless?} 
\markright{The Case of the Second Smallest Modulus}
\author{Katalin Gyarmati}

\maketitle

\begin{abstract}
 We present a concise, elementary proof that the second smallest modulus of a covering system must be bounded, simplifying the known  results through a maximality argument.
\end{abstract}

{\let\thefootnote\relax}

\section{The Journey of a \$10 Conjecture}

The theory of covering congruences was one of Paul Erdős's most beloved research areas. The story began in 1950 when Erdős \cite{erdos} studied a problem posed by Romanov: does there exist a set of positive density consisting of integers that \textit{cannot} be represented in the form $p + 2^k$? Erdős's answer was yes, and for the proof, he used a construction where every integer was covered by a finite number of congruence classes whose moduli were distinct and greater \mbox{than 1.}

\begin{dfn}
A finite set of congruences $\mathcal{M} = \{a_i \pmod{m_i}\}_{i=1}^k$ is called a \textit{covering system} if every integer satisfies at least one of the congruences. The system is \textit{non-trivial} if the moduli are greater than one and distinct ($1 < m_1 < m_2 < \dots < m_k$).
\end{dfn}

It was then that the question arose which has intrigued number theorists for decades: does there exist a covering system with an arbitrarily large smallest modulus ($m_1$)? Erdős was so confident in the answer that he initially offered 10, and later 1000 dollars for the solution. He was convinced that $m_1$ could be as large as desired.

For decades, only modest progress was made. In the early 2000s, probabilistic methods, especially the Lovász Local Lemma, brought new hope. The essence of the lemma is that if the ``bad'' events (for example, an integer remaining uncovered) are only weakly dependent, then there is a positive chance for a covering. However, the community had to wait until 2015 for complete success.

Bob Hough \cite{hough} shook the foundations of the field when he disproved Erdős's conjecture: $m_1$ cannot grow without bound. 
Hough's iterative sieving procedure showed that in every covering system, $m_1$ is bounded ($m_1 \le 10^{16}$). This bound was narrowed down to the currently best known $616,000$ in 2022 by Balister, Bollobás, Morris, Sahasrabudhe, and Tiba \cite{balister}, using a refinement of the Lovász Local Lemma.

But what happens to the next modulus in line? Can a covering system keep the value of $m_1$ low (for example, $m_1=2$) while pushing the value of $m_2$ to infinity? Although the literature \cite{klein, cummings} has already answered this question for the general $j$-th modulus, the methods applied are often quite complex.

In this note, we provide an elementary, concise one-page proof for the boundedness of $m_2$. Our approach avoids complex sieve-theoretic machinery, building instead on the internal contradictions of the structure and the self-similarity of covering systems. We show that if we accept the boundedness of $m_1$, then the ``escape'' of $m_2$ is a mathematical impossibility.

\section{The Proof}

We assume, for the sake of contradiction, that the second smallest modulus, $m_2$, can be arbitrarily large. 

\bigskip
Denote by $C_1$ the known upper bound for the first modulus (e.g., $616,000$). We consider the set of all covering systems $\mathcal{M} = \{a_i \pmod{m_i}\}_{i=1}^k$ such that $m_2 > C_1^2$. Since $m_1 \le C_1$, there are only finitely manypossible values for $m_1$. Among all such systems, we choose one in which $m_1$ is maximal. Furthermore, we assume that the moduli are ordered by magnitude such that $m_1<m_2<m_3<\dots<m_k$ holds.

\bigskip
By our indirect assumption, there exists a covering system
$$\mathcal{T} = \{b_s \pmod{T_s}\}_{s=1}^v,$$
where the second modulus $T_2$ is arbitrarily large. We choose $T_2$ such that:
$$T_2 > \frac{m_k}{m_1},$$
where $m_k$ is the largest modulus of our original system $\mathcal{M}$.

\bigskip
Replace the first congruence of $\mathcal{M}$, $a_1 \pmod{m_1}$, with a ``spread'' version of the system $\mathcal{T}$:
$$\{a_1 + b_s m_1 \pmod{m_1 T_s} : 1 \le s \le v\}.$$
The new system, denoted by $\mathcal{M}'$, remains a covering system. Its set of moduli is:
$$\mathcal{M}_{\text{mod}}= \{m_2, m_3, \dots, m_k\} \cup \{m_1 T_1, m_1 T_2, \dots, m_1 T_v\}.$$

\bigskip 
The new moduli are clearly distinct. Since the moduli $T_s$ are distinct, the values $m_1 T_s$ are also distinct. Furthermore, since $T_1 > 1$, we have $m_1 T_1 > m_1$. We now examine the ordering of the new moduli:
\begin{enumerate}
    \item $m_1 T_1 < m_1 T_2$, and by our choice, $m_1 T_2 > m_k$.
    \item Since $m_1 \le C_1$ and $T_1$ is the first modulus of a covering system (thus $T_1 \le C_1$), we have $m_1 T_1 \le C_1^2 < m_2$. Thus, the order of the moduli is
        $$m_1T_1<m_2<m_3<\dots<m_k<m_1T_2<m_1T_3<\dots<m_1T_v.$$
\end{enumerate}
This implies that in the new system, the minimum modulus is $m_1' = m_1 T_1$. Since $m_2' =m_2 > C_1^2$, the new system $\mathcal{M}'$ belongs to the same class of covering systems. However, $m_1' = m_1 T_1 > m_1$ directly contradicts the maximality of $m_1$ in the original system. \qed

\section*{Funding Information}
Research supported by the Hungarian National Research Development and 
Innovation Fund KKP133819.

\bigskip
\noindent Katalin Gyarmati \\
Eötvös Loránd University, Budapest, Hungary,\\ katalin.gyarmati@gmail.com

\end{document}